\documentclass[twosided,reqno]{amsart}

\usepackage{amsmath}
\usepackage{amssymb}
\usepackage{amsthm}

\theoremstyle{theorem}
\newtheorem{theorem}{\scshape Theorem }[section]
\newtheorem{lemma}[theorem]{\scshape Lemma}
\newtheorem{corollary}[theorem]{\scshape Corollary}

\theoremstyle{definition}

\numberwithin{equation}{section}

\begin{document}

\title[Associated sequences of special polynomials]{Some identities involving associated sequences of special polynomials}

\author{Taekyun Kim$^1$}
\address{$^1$ Department of Mathematics, Kwangwoon University, Seoul 139-701, Republic of Korea.}
\email{tkkim@kw.ac.kr}

\author{Dae San Kim$^2$}
\address{$^2$ Department of Mathematics, Sogang University, Seoul 121-742, Republic of Korea.}
\email{dskim@sogang.ac.kr}


\maketitle

\begin{abstract}
In this paper, we study some properties of associated sequences of special polynomials. From the properties of associated sequences of polynomials, we derive some interesting identities of special polynomials.
\end{abstract}

\section{Introduction}

For $r \in {\mathbb{R}}$,
the {\it{Bernoulli polynomials}} of order $r$ are defined by the generating function to be
\begin{equation}\label{1}
\left(\frac{t}{e^t-1}\right)^re^{xt}=\sum_{n=0} ^{\infty} B_n ^{(r)} (x) \frac{t^n}{n!},~(r \in {\mathbb{R}}),{\text{ (see [12,13,14,18,21])}}.
\end{equation}
In the special case, $x=0$, $B_n ^{(r)} (0)=B_n ^{(r)}$ are called the $n$-th Bernoulli numbers of order $r$. It is also well known that the {\it{Euler polynomials}} of order $r$ are defined by the generating function to be
\begin{equation}\label{2}
\left(\frac{2}{e^t+1}\right)^re^{xt}=\sum_{n=0} ^{\infty}E_n ^{(r)} (x)\frac{t^n}{n!},~(r \in {\mathbb{R}}) ,{\text  { \text  (see [9,10,11,19,20])}}.
\end{equation}
Let $x=0$. Then $E_n ^{(r)}(0)=E_n ^{(r)}$ are called the {\it{$n$-th Euler numbers}} of order $r$.

Let ${\mathcal{F}}$ be the set of all formal power series in the variable $t$ over ${\mathbb{C}}$ with
\begin{equation}\label{3}
{\mathcal{F}}=\left\{ \left.f(t)=\sum_{k=0} ^{\infty} \frac{a_k}{k!} t^k~\right|~ a_k \in {\mathbb{C}} \right\}.
\end{equation}
Let ${\mathbb{P}}$ be the algebra of polynomials in the variable $x$ over ${\mathbb{C}}$ and ${\mathbb{P}}^{*}$ be the vector space  of all linear functionals on ${\mathbb{P}}$. The action of the linear functional $L$ on a polynomial $p(x)$ is defined by $\left< L~|~p(x)\right>$ and the vector space structure on ${\mathbb{P}}^{*}$ is derived by $\left<L+M|p(x) \right>=\left<L|p(x) \right>+\left<M|p(x) \right>$, $\left<cL|p(x) \right>=c\left<L|p(x) \right>$, where $c$ is a complex constant.

For  $f(t)=\sum_{k=0} ^{\infty} \frac{a_k}{k!}t^k \in {\mathcal{F}}$, we define a linear functional on ${\mathbb{P}}$ by setting
\begin{equation}\label{4}
\left<f(t)|x^n \right>=a_n,~ (n\geq 0),{\text{ (see [3,8,17])}}.
\end{equation}
By \eqref{3} and \eqref{4}, we get
\begin{equation}\label{5}
\left<t^k | x^n \right>=n! \delta_{n,k},~(n,k \geq 0),{\text{ (see [4,5,7,10,17,18])}},
\end{equation}
where $\delta_{n,k}$ is the Kronecker symbol.

Let $f_L(t)=\sum_{k=0} ^{\infty} \frac{\left<L|x^k\right>}{k!}t^k$. Then, by \eqref{5}, we get $\left<f_L(t)|x^n\right>=\left<L|x^n\right>$. So, we see that $f_L(xt)=L$. The map $L \mapsto f_L(t)$ is a vector space isomorphism from ${\mathbb{P}}^{*}$ onto ${\mathcal{F}}$. Henceforth, ${\mathcal{F}}$ is thought of as both a formal power series and a linear functional (see \cite{04,06,08,17,18}). We call ${\mathcal{F}}$ the {\it{umbral algebra}}. The umbral calculus is the study of umbral algebra (see \cite{04,10,17}).

The order $o(f(t))$ of the non-zero power series $f(t)$ is the smallest integer $k$ for which the coefficient of $t^k$ does not vanish. If $o(f(t))=1$, then $f(t)$ is called a {\it{delta series}} and if $o(f(t))=0$, then $f(t)$ is called an {\it{invertible series}}. Let $o(f(t))=1$ and $o(g(t))=0$. Then there  exists a unique sequence $S_n(x)$ of polynomials such that $\left<g(t)f(t)^k|S_n(x)\right>=n!\delta_{n,k}$ where $n,k \geq 0$. The sequence $S_n(x)$ is called {\it{Sheffer sequence}} for $(g(t),f(t))$, which is denoted by $S_n(x)\sim (g(t),f(t))$. If $S_n(x)\sim(1,f(t))$, then $S_n(x)$ is called the {\it{associated sequence}} for $f(t)$. By \eqref{5}, we see that $\left. \left<e^{yt} \right| p(x)\right>=p(y)$.

Let  $f(t)\in{\mathcal{F}}$ and $p(x)\in{\mathbb{P}}$. Then we have
\begin{equation}\label{6}
f(t)=\sum_{k=0} ^{\infty} \frac{\left<f(t)|x^k\right>}{k!}t^k,~p(x)=\sum_{k=0} ^{\infty} \frac{\left<t^k|p(x)\right>}{k!}x^k, {\text{ (see [10,17])}},
\end{equation}
and
\begin{equation}\label{7}
\left. \left<f_1(t)f_2(t)\cdots f_m(t)\right| x^n \right>=\sum_{i_1+\cdots +i_m=n} \binom{n}{i_1,\ldots,i_m}\left.\left<f_1(t)\right|x^{i_1}\right>\cdots \left.\left<f_m(t)\right|x^{i_1}\right>,
\end{equation}
where $f_1(t),f_2(t),\cdots ,f_m(t)\in {\mathcal{F}}$,  (see \cite{04,10,17}).
From \eqref{6}, we have
\begin{equation}\label{8}
p^{(k)}(0)=\left<t^k|p(x)\right>,~\left< 1\left|p^{(k)}(x)\right.\right>=p^{(k)}(0).
\end{equation}
Thus, by \eqref{8}, we get
\begin{equation}\label{9}
t^k p(x)=p^{(k)}(x)=\frac{d^kp(x)}{dx^k},~(k \geq 0),{\text{ (see [17])}}.
\end{equation}
For $S_n(x)\sim(g(t),f(t))$, we have the following equations:
\begin{equation}\label{10}
S_n(x+y)=\sum_{k=0} ^n \binom{n}{k} p_k(y)S_{n-k} (x),{\text{ where }}p_k(y)=g(t)S_k(y),
\end{equation}
and
\begin{equation}\label{11}
\frac{1}{g({\bar{f}}(t))}e^{y{\bar{f}}(t)}=\sum_{k=0} ^{\infty} \frac{S_k(y)}{k!}t^k,{\text{ for all }}y \in {\mathbb{C}}{\text{ (see [10,15,16,17,18])}},
\end{equation}
where ${\bar{f}}(t)$ is the compositional inverse of $f(t)$.

Let $p_n(x)\sim(1,f(t))$ and $q_n(x)\sim(1,g(t))$. Then the transfer formula for associated sequence implies that, for $n\in\Bbb N$,
\begin{equation}\label{12}
q_n(x)=x\left(\frac{f(t)}{g(t)}\right)^nx^{-1}p_n(x),{\text{ (see [11,17,22])}}.
\end{equation}
Now we introduce several important sequences which are used to derive our results in this paper (see \cite{10, 11,17}):

(The Poisson-Charlier sequences)

\begin{equation*}
C_n(x;a)=\sum_{k=0} ^n \binom{n}{k}(-1)^{n-k}a^{-k}(x)_k\sim\left(e^{a(e^t-1)},a(e^t-1)\right),
\end{equation*}
where $a \neq 0$, $(x)_n=x(x-1)\cdots(x-n+1)$,
\begin{equation}\label{13}
\sum_{k=0}^nC_n(k;a)\frac{t^k}{k!}e^{-t}=\left(\frac{t-a}{a}\right)^n,~(a \neq 0),~n \in {\mathbb{N}\cup \{0\}},
\end{equation}
(The Abel sequences)

\begin{equation}\label{14}
A_n(x;b)=x(x-bn)^{n-1}\sim\left(1,te^{bt}\right),~(b \neq 0),
\end{equation}
(The Mittag-Leffler sequences)

\begin{equation}\label{15}
M_n(x)=\sum_{k=0} ^n \binom{n}{k} (n-1)_{n-k} 2^k (x)_k \sim \left(1,\frac{e^t-1}{e^t+1}\right),
\end{equation}
(The exponential sequences)
\begin{equation}\label{16}
\phi_n(x)=\sum_{k=0} ^n S_2(n,k)x^k \sim\left(1, \log (1+t)\right),
\end{equation}
and

(The Laguerre sequences)

\begin{equation}\label{17}
L_n(x)=\sum_{k=1} ^n \binom{n-1}{k-1}\frac{n!}{k!}(-x)^k\sim\left(1,\frac{t}{t-1}\right).
\end{equation}

In this paper, we study some properties of associated sequences of special polynomials. From the properties of associated sequences of specials polynomials, we derive some interesting identities involving associated sequences of special polynomials.

\section{Associated sequences of special polynomials.}

As is well known, the Bessel differential equation is given by
\begin{equation}\label{18}
x^2 y^{''}+2(x+1)y^{'}+n(n+1)y=0,{\text{ (see [1,2])}}.
\end{equation}
From \eqref{18}, we have the solution of \eqref{18} as follows:
\begin{equation}\label{19}
y_n (x)=\sum_{k=0} ^n \frac{(n+k)!}{(n-k)!k!}\left(\frac{x}{2}\right)^k,{\text{ (see [1,2])}}.
\end{equation}
Let us consider the following associated sequences:
\begin{equation}\label{20}
p_n(x)\sim\left(1,t-\frac{t^2}{2}\right),~x^n \sim (1,t),{\text{ (see [1,2,10,17])}}.
\end{equation}
From \eqref{12} and \eqref{20}, for $n\in \Bbb N$, we have
\begin{equation}\label{21}
\begin{split}
p_n(x)&=x\left(\frac{t}{t-\frac{t^2}{2}}\right)^nx^{-1}x^n=x\left(1-\frac{t}{2}\right)^{-n}x^{n-1} \\
&=x\sum_{k=0} ^{\infty}\binom{-n}{k}(-1)^k\left(\frac{t}{2}\right)^kx^{n-1} \\
&=x\sum_{k=0} ^{n-1} \binom{n+k-1}{k} \left(\frac{1}{2}\right)^k(n-1)_kx^{n-1-k} \\
&=\sum_{k=0} ^{n-1} \frac{(n+k-1)!}{k!(n-1-k)!}\left(\frac{1}{2}\right)^kx^{n-k}\\
&=\sum_{k=1} ^n \frac{(2n-k-1)!}{(n-k)!(k-1)!}\left(\frac{1}{2}\right)^{n-k}x^k.
\end{split}
\end{equation}
By \eqref{19} and \eqref{21}, we get
\begin{equation}\label{22}
p_n(x)=x^n y_{n-1} \left(\frac{1}{x}\right)\sim\left(1,t-\frac{t^2}{2}\right).
\end{equation}
From \eqref{11} and \eqref{22}, we can derive the following generating function of $p_n(x)$:
\begin{equation}\label{23}
\sum_{k=0} ^{\infty}p_k(x)\frac{t^k}{k!}=e^{x\left(1-(1-2t)^{\frac{1}{2}}\right)},
\end{equation}
and, by \eqref{10}, we get
\begin{equation}\label{24}
(x+y)^ny_{n-1}\left(\frac{1}{x+y}\right)=\sum_{k=0} ^n \binom{n}{k}x^ky^{n-k}y_{k-1}\left(\frac{1}{x}\right)y_{n-k-1}\left(\frac{1}{y}\right).
\end{equation}
By \eqref{12} and \eqref{20}, we get
\begin{equation}\label{25}
x^n=x\left(\frac{t-\frac{t^2}{2}}{t}\right)^nx^{-1} p_n(x)=x\left(\frac{t-2}{-2}\right)^nx^{-1}p_n(x).
\end{equation}
Thus, by \eqref{13}, \eqref{21} and \eqref{25}, we get
\begin{equation}\label{26}
\begin{split}
(-1)^nx^{n-1}&=\left(\frac{t-2}{2}\right)^2x^{-1}p_n(x)=\sum_{k=0} ^{\infty}C_n(k;2)\frac{t^k}{k!}e^{-t}\left(x^{-1}p_n(x)\right) \\
&=\sum_{k=0} ^{n-1} C_n(k;2)\frac{t^k}{k!}(x-1)^{-1}p_n(x-1)\\
&=\sum_{k=0} ^{n-1} C_n(k;2)\frac{t^k}{k!}\sum_{l=1} ^n \frac{(2n-l-1)!}{(l-1)!(n-l)!}\left(\frac{1}{2}\right)^{n-l}(x-1)^{l-1}\\
&=\sum_{k=0} ^{n-1} C_n(k;2)\sum_{l=k+1} ^n \frac{(2n-l-1)!}{(l-1)!(n-l)!}\binom{l-1}{k}\left(\frac{1}{2}\right)^{n-l}(x-1)^{l-1-k}\\
&=\sum_{m=0} ^{n-1}\sum_{k=0} ^{n-m-1}C_n(k;2)\binom{m+k}{k}\frac{(2n-m-k-2)!}{(m+k)!(n-m-k-1)!}\left(\frac{1}{2}\right)^{n-m-k-1}(x-1)^m.
\end{split}
\end{equation}
From \eqref{25}, we have
\begin{equation}\label{27}
\begin{split}
x^{n-1}&=\left(1-\frac{t}{2}\right)^nx^{-1}p_n(x)=\sum_{k=0} ^n \binom{n}{k}\left(-\frac{t}{2}\right)^kx^{-1}p_n(x)\\
&=\sum_{k=0} ^{n-1}\sum_{l=k+1} ^n \binom{n}{k}(l-1)_k\frac{(2n-l-1)!}{(l-1)!(n-l)!}(-1)^k\left(\frac{1}{2}\right)^{n+k-l}x^{l-1-k} \\
&=\sum_{m=0} ^{n-1} \sum_{k=0} ^{n-m-1} \binom{n}{k}(m+k)_k \frac{(2n-m-k-2)!(-1)^k}{(m+k)!(n-m-1-k)!}\left(\frac{1}{2}\right)^{n-m-1}x^m\\
&=\sum_{m=0} ^{n-1}\left\{\sum_{k=0} ^{n-m-1}(-1)^k \left(\frac{1}{2}\right)^{n-m-1} \binom{n}{k}\frac{(2n-m-k-2)!}{m!(n-m-k-1)!}\right\}x^m.
\end{split}
\end{equation}
Therefore, by \eqref{26} and \eqref{27}, we obtain the following theorem.
\begin{theorem}\label{thm1}
For $n \in {\mathbb{N}}$, we have
\begin{equation*}
(-1)^nx^{n-1}=\sum_{m=0} ^{n-1}\sum_{k=0} ^{n-m-1}C_n(k;2)\binom{m+k}{k}\frac{(2n-m-k-2)!}{(m+k)!(n-m-k-1)!}\left(\frac{1}{2}\right)^{n-m-k-1}(x-1)^m.
\end{equation*}
Moreover,
\begin{equation*}
\sum_{k=0} ^{n-m-1}(-1)^k\left(\frac{1}{2}\right)^{n-m-1} \binom{n}{k}\frac{(2n-m-k-2)!}{m!(n-m-k-1)!}=0,
\end{equation*}
where $0 \leq m \leq n-2$.
\end{theorem}
Let us consider the following associated sequences:
\begin{equation}\label{28}
p_n(x)\sim\left(1,te^{c(e^t-1)}\right),~c \neq 0,~A_n(x;b)=x(x-bn)^{n-1}\sim\left(1,te^{bt}\right),~b \neq 0.
\end{equation}
By \eqref{12} and \eqref{28}, we get
\begin{equation}\label{29}
p_n(x)=x \left(\frac{t}{te^{c(e^t-1)}}\right)^nx^{-1}x^n=x\sum_{k=0} ^{\infty}\frac{(-nc)^k}{k!}(e^t-1)^kx^{n-1}.
\end{equation}
We recall that Newton's difference operator $\Delta$ is defined by $\Delta f(x)=f(x+1)-f(x)$. For $n \in {\mathbb{N}}$, we easily see that
\begin{equation}\label{30}
\Delta^np(x)=\sum_{k=0} ^n \binom{n}{k}(-1)^{n-k} p(x+k).
\end{equation}
By \eqref{30}, we get
\begin{equation}\label{31}
(e^t-1)^kp(x)=\sum_{l=0} ^k \binom{k}{l}(-1)^{k-l}e^{lt}p(x)=\sum_{l=0} ^k \binom{k}{l}(-1)^{l-k}p(x+l)=\Delta^kp(x).
\end{equation}
In particular, if we take $p(x)=x^{n-1}$, then we have
\begin{equation}\label{32}
(e^t-1)^kx^{n-1}=\sum_{j=0} ^k \binom{k}{j} (-1)^{k-j} (x+j)^{n-1}.
\end{equation}
From \eqref{29} and \eqref{32}, we have
\begin{equation}\label{33}
p_n(x)=x\sum_{k=0} ^{n-1}\sum_{j=0} ^k \frac{(-1)^j(nc)^k}{k!}\binom{k}{j}(x+j)^{n-1}\sim\left(1,te^{c(e^t-1)}\right), ~c \neq 0.
\end{equation}
Therefore, by \eqref{33}, we obtain the following lemma.
\begin{lemma}\label{lemma2}
For $c \neq 0$ and $n \in {\mathbb{N}}$, let $p_n(x)\sim\left(1,te^{c(e^t-1)}\right)$. Then we have
\begin{equation*}
p_n(x)=x\sum_{k=0} ^{n-1} \sum_{j=0} ^k \frac{(-1)^j(nc)^k}{k!}\binom{k}{j}(x+j)^{n-1}.
\end{equation*}
\end{lemma}
From the definition of Abel sequences and \eqref{28}, we note that
\begin{equation}\label{34}
\begin{split}
A_n(x;b)&=x(x-bn)^{n-1}=x\left(\frac{te^{c(e^t-1)}}{te^{bt}}\right)^nx^{-1}p_n(x) \\
&=x\left(e^{-bt-c(1-e^t)}\right)^nx^{-1}p_n(x)=x\left(\sum_{l=0} ^{\infty}\frac{a_l ^{(-b)}(-c)}{l!}t^l \right)^nx^{-1}p_n(x),
\end{split}
\end{equation}
where $a_n ^{(\beta)}(x)\sim \left((1-t)^{-\beta},\log (1-t)\right)$ is the {\it{actuarial polynomial}} with the generating function  given by
\begin{equation*}
\sum_{l=0} ^{\infty}\frac{a_l ^{(\beta)}(x)}{l!}t^l=e^{\beta t+x(1-e^t)}.
\end{equation*}
By Lemma \ref{lemma2} and \eqref{34}, we get
\begin{equation}\label{35}
\begin{split}
&A_n(x;b) \\
=&x \left\{\sum_{m=0} ^{\infty} \sum_{l_1+\cdots+l_n=m}\binom{m}{l_1,\ldots,l_n}a_{l_1} ^{(-b)}(-c)\cdots a_{l_n} ^{(-b)}(-c) \frac{t^m}{m!}\right\}\\
& \times \left\{\sum_{k=0} ^{n-1}\sum_{j=0} ^k \frac{(-1)^j(nc)^k}{k!}\binom{k}{j}(x+j)^{n-1}\right\} \\
=&x\sum_{m=0} ^{n-1}\sum_{l_1+\cdots+l_n=m}\sum_{k=0} ^{n-1}\sum_{j=0} ^k\binom{n-1}{m}\binom{m}{l_1,\ldots,l_n}\binom{k}{j}\left(\prod_{j=1} ^n a_{l_j} ^{(-b)}(-c)\right)\left(\frac{(-1)^j(nc)^k}{k!}\right)(x+j)^{n-1-m}.
\end{split}
\end{equation}
Therefore, by \eqref{35}, we obtain the following theorem.
\begin{theorem}\label{thm3}
For $n \geq 1$, $b \neq 0$, $c \neq 0$, we have
\begin{equation*}
\begin{split}
&A_n(x;b)\\
=&x\sum_{m=0} ^{n-1}\sum_{l_1+\cdots+l_n=m}\sum_{k=0} ^{n-1}\sum_{j=0} ^k\binom{n-1}{m}\binom{m}{l_1,\ldots,l_n}\binom{k}{j}\left(\prod_{j=1} ^n a_{l_j} ^{(-b)}(-c)\right)\left(\frac{(-1)^j(nc)^k}{k!}\right)(x+j)^{n-1-m}.
\end{split}
\end{equation*}
\end{theorem}
For \eqref{34}, we note that
\begin{equation}\label{36}
\begin{split}
A_n(x;b)&=x\left(e^{c(e^t-1)}\right)^ne^{-nbt}x^{-1}p_n(x)\\
&=x\left(e^{c(e^t-1)}\right)^n(x-nb)^{-1}p_n(x-nb).
\end{split}
\end{equation}
By \eqref{16} and Lemma \ref{lemma2}, we easily see that the generating function of exponential sequences is given by
\begin{equation}\label{37}
\sum_{k=0} ^{\infty} \phi_k(x)\frac{t^k}{k!}=e^{x(e^t-1)}.
\end{equation}
From \eqref{36} and \eqref{37}, we have
\begin{equation}\label{38}
\begin{split}
&A_n(x;b) \\
&=x\left\{\sum_{m=0} ^{\infty}\sum_{l_1+\cdots+l_n=m}\binom{m}{l_1,\ldots,l_n}\left(\prod_{j=1} ^n \phi_{l_j} (c)\right)\frac{t^m}{m!}\right\}(x-nb)^{-1}p_n(x-nb) \\
&=x\sum_{m=0} ^{\infty}\sum_{l_1+\cdots+l_n=m}\binom{m}{l_1,\ldots,l_n}\left(\prod_{j=1} ^n \phi_{l_j} (c)\right)\frac{t^m}{m!}\sum_{k=0} ^{n-1} \sum_{j=0} ^k \frac{(-1)^j(nc)^k}{k!}
\binom{k}{j}(x-nb+j)^{n-1} \\
&=x\sum_{m=0} ^{n-1}\sum_{l_1+\cdots+l_n=m}\sum_{k=0} ^{n-1} \sum_{j=0} ^k \binom{n-1}{m}\binom{m}{l_1,\ldots,l_n}\binom{k}{j}\left(\prod_{j=1} ^n \phi_{l_j} (c)\right)\frac{(-1)^j(nc)^k}{k!}(x-nb+j)^{n-m-1}.
\end{split}
\end{equation}
Therefore, by \eqref{38}, we obtain the following corollary.
\begin{corollary}\label{coro4}
For $n \geq 1$, $b \neq 0$, $c \neq 0$, we have
\begin{equation*}
\begin{split}
&A_n(x;b) \\
=&x\sum_{m=0} ^{n-1}\sum_{l_1+\cdots+l_n=m}\sum_{k=0} ^{n-1} \sum_{j=0} ^k \binom{n-1}{m}\binom{m}{l_1,\ldots,l_n}\binom{k}{j}\left(\prod_{j=1} ^n \phi_{l_j} (c)\right)\frac{(-1)^j(nc)^k}{k!}(x-nb+j)^{n-m-1}.
\end{split}
\end{equation*}
\end{corollary}
Note that $x^n \sim (1,t)$. By \eqref{12}, \eqref{13} and \eqref{17}, we get
\begin{equation}\label{39}
\begin{split}
L_n(x)&=x\left(\frac{t}{\frac{t}{t-1}}\right)^nx^{-1}x^n=x(t-1)^nx^{n-1} \\
&=x\left(\sum_{k=0} ^{n-1}C_n(k;1)\frac{t^k}{k!}e^{-t}\right)x^{n-1}=x\sum_{k=0} ^{n-1}C_n(k;1)\frac{t^k}{k!}(x-1)^{n-1}\\
&=x\sum_{k=0} ^{n-1}C_n(k;1) \binom{n-1}{k}(x-1)^{n-1-k}=x\sum_{k=0} ^{n-1}\binom{n-1}{k}C_n(n-1-k;1)(x-1)^k.
\end{split}
\end{equation}
Therefore, by \eqref{39}, we obtain the following theorem.
\begin{theorem}\label{thm5}
For $n \geq 1$, we have
\begin{equation*}
L_n(x)=x\sum_{k=0} ^{n-1}\binom{n-1}{k}C_n(n-1-k;1)(x-1)^k.
\end{equation*}
\end{theorem}
 Mott considered the associated sequences for $f(t)=\frac{-2t}{1-t^2}$. That is, the Mott sequence is given by
\begin{equation}\label{40}
S_n(x)\sim\left(1,\frac{-2t}{1-t^2}\right).
\end{equation}
From \eqref{40}, we note that the generating function of Mott sequences is given by
\begin{equation*}
\sum_{k=0} ^{\infty}S_k(x)\frac{t^k}{k!}=\exp\left(x\left(\frac{1-\sqrt{1+t^2}}{t}\right)\right).
\end{equation*}
By \eqref{12}, \eqref{17} and \eqref{40}, we get
\begin{equation}\label{41}
\begin{split}
S_n(x)&=x\left(\frac{\frac{t}{t+1}}{\frac{-2t}{1-t^2}}\right)^nx^{-1}L_n(-x)=x\left(\frac{t-1}{2}\right)^n x^{-1}L_n(-x) \\
&=2^{-n}x(t-1)^nx^{-1}L_n(-x)=2^{-n}x\left(\sum_{k=0} ^{n-1}C_n(k;1)\frac{t^k}{k!}e^{-t}\right)x^{-1}L_n(-x) \\
&=2^{-n}x\sum_{k=0} ^{n-1}C_n(k;1)\frac{t^k}{k!}(x-1)^{-1}L_n(1-x) \\
&=2^{-n}x\sum_{k=0} ^{n-1}C_n(k;1)\frac{1}{k!}\sum_{l=1} ^n \binom{n-1}{l-1}\frac{n!}{l!}t^k(x-1)^{l-1} \\
&=\frac{n!}{2^n}\sum_{k=0} ^{n-1} \sum_{l=1} ^n \binom{n-1}{l-1}\binom{l-1}{k}\frac{C_n(k;1)}{l!}x(x-1)^{l-1-k}.
\end{split}
\end{equation}
Thus, by \eqref{41}, we obtain the following lemma.
\begin{lemma}\label{lemma6}
For $n \in {\mathbb{N}}$, let $S_n(x)\sim\left(1,\frac{-2t}{1-t^2}\right)$. Then we have
\begin{equation*}
S_n(x)=\frac{n!}{2^n}\sum_{k=0} ^{n-1} \sum_{l=1} ^n \binom{n-1}{l-1}\binom{l-1}{k}\frac{C_n(k;1)}{l!}x(x-1)^{l-1-k}.
\end{equation*}
\end{lemma}
As is known, we have
\begin{equation}\label{42}
xB_{n-1} ^{(an)}(x)\sim \left(1,t\left(\frac{e^t-1}{t}\right)^a\right),~xE_{n-1} ^{(bn)}(x)\sim\left(1,t\left(\frac{e^t+1}{2}\right)^b \right),
\end{equation}
where $a,b$ are positive integers (see \cite{10,11,17}). For $n \geq 1$, by \eqref{12} and \eqref{42}, we get
\begin{equation}\label{43}
\begin{split}
xE_{n-1} ^{(bn)}(x)&=x\left(\frac{t\left(\frac{e^t-1}{t}\right)^a}{t\left(\frac{e^t+1}{2}\right)^b}\right)^nx^{-1}B_{n-1} ^{(an)}(x) \\
&=x\left(\frac{e^t+1}{2}\right)^{-bn}\left(\frac{e^t-1}{t}\right)^{an}B_{n-1} ^{(an)}(x).
\end{split}
\end{equation}
Thus, by \eqref{43}, we get
\begin{equation}\label{44}
\left(\frac{e^t+1}{2}\right)^{bn}E_{n-1} ^{(bn)} (x)=\left(\frac{e^t-1}{t}\right)^{an}B_{n-1} ^{(an)}(x).
\end{equation}
\begin{equation}\label{45}
\begin{split}
{\text{LHS of \eqref{44}}}&= 2^{-bn}\left(e^{t}+1\right)^{bn}E_{n-1} ^{(bn)}(x)=2^{-bn}\sum_{k=0} ^{bn}\binom{bn}{k}e^{kt}E_{n-1} ^{(bn)}(x) \\
&=2^{-bn}\sum_{k=0} ^{bn}\binom{bn}{k}E_{n-1} ^{(bn)} (x+k).
\end{split}
\end{equation}
\begin{equation}\label{46}
\begin{split}
{\text{RHS of \eqref{44}}}&= \left(\frac{1}{t}\right)^{an}(an)!\sum_{l=an} ^{\infty}S_2(l,an)\frac{t^l}{l!}B_{n-1} ^{(an)}(x) \\
&=\sum_{l=0} ^{n-1} \frac{(an)!}{(l+an)!}S_2(l+an,an)(n-1)_lB_{n-1-l} ^{(an)}(x) \\
&=(n-1)!\sum_{l=0} ^{n-1}\frac{(an)!}{(l+an)!(n-1-l)!}S_2(l+an,an)B_{n-1-l} ^{(an)}(x),
\end{split}
\end{equation}
where $S_2(n,k)$ is the Stirling number of the second kind. Therefore, by \eqref{45} and \eqref{46}, we obtain the following theorem.
\begin{theorem}\label{thm7}
For $n \geq 1$, $a,b \in {\mathbb{N}}\cup\left\{0\right\}$, we have
\begin{equation*}
\sum_{k=0} ^{bn} \binom{bn}{k}E_{n-1} ^{(bn)}(x+k)=2^{bn}(n-1)!\sum_{l=0} ^{n-1} \frac{(an)!}{(l+an)!(n-1-l)!}S_2(l+an,an)B_{n-1-l} ^{(an)}(x).
\end{equation*}
\end{theorem}
The Pidduck sequences is given by
\begin{equation}\label{47}
P_n(x)\sim \left(\frac{2}{e^t+1},\frac{e^t-1}{e^t+1}\right).
\end{equation}
From \eqref{47}, we can derive the generating function of the Pidduck sequences as follows:
\begin{equation}\label{48}
\sum_{k=0} ^{\infty}P_k(x)\frac{t^k}{k!}=(1-t)^{-1}\left(\frac{1+t}{1-t}\right)^x.
\end{equation}
Let $S_n(x)\sim \left(1,\frac{2t}{e^t+1}\right)$. Then, from \eqref{12}, \eqref{15} and \eqref{47}, we have
\begin{equation}\label{49}
\begin{split}
M_n(x)&=\frac{2}{e^t+1}P_n(x)=x\left(\frac{\frac{2t}{e^t+1}}{\frac{e^t-1}{e^t+1}}\right)^nx^{-1}S_n(x) \\
&=2^nx\left(\frac{t}{e^t-1}\right)^nx^{-1}S_n(x).
\end{split}
\end{equation}
By \eqref{12}, we easily get
\begin{equation}\label{50}
\begin{split}
S_n(x)&=x\left(\frac{\frac{2t}{e^t+1}}{t}\right)^{-n}x^{-1}x^n=2^{-n}x(e^t+1)^nx^{n-1} \\
&=2^{-n}x\sum_{j=0} ^n \binom{n}{j}e^{jt}x^{n-1}=2^{-n}x\sum_{j=0} ^n\binom{n}{j}(x+j)^{n-1}.
\end{split}
\end{equation}
From \eqref{49} and \eqref{50}, we have
\begin{equation}\label{51}
\begin{split}
M_n(x)&=2^nx\left(\frac{t}{e^t-1}\right)^n\left(2^{-n}\sum_{j=0} ^n \binom{n}{j}(x+j)^{n-1}\right) \\
&=\sum_{j=0} ^n \binom{n}{j}x\left(\frac{t}{e^t-1}\right)^n(x+j)^{n-1}=\sum_{j=0} ^n \binom{n}{j}xB_{n-1} ^{(n)}(x+j).
\end{split}
\end{equation}
By \eqref{49} and \eqref{51}, we get
\begin{equation}\label{52}
\begin{split}
P_n(x)&=\frac{1}{2}(e^t+1)M_n(x)=\frac{1}{2}(e^t+1)\sum_{j=0} ^n \binom{n}{j}xB_{n-1} ^{(n)}(x+j) \\
&=\frac{1}{2}\sum_{j=0} ^n \binom{n}{j}\left\{(x+1)B_{n-1} ^{(n)}(x+1+j)+xB_{n-1} ^{(n)}(x+j)\right\} \\
&=\frac{1}{2}\left\{\sum_{j=1} ^n\left(\binom{n}{j-1}(x+1)+\binom{n}{j}x\right)B_{n-1} ^{(n)}(x+j)+(x+1)B_{n-1} ^{(n)}(x+n+1)+xB_{n-1} ^{(n)}(x)\right\} \\
&=\frac{1}{2}\sum_{j=0} ^{n+1}\left\{\binom{n+1}{j}x+\binom{n}{j-1}\right\}B_{n-1} ^{(n)}(x+j).
\end{split}
\end{equation}
Therefore, by \eqref{52}, we obtain the following theorem.
\begin{theorem}\label{thm8}
For $n \geq 1$, we have
\begin{equation*}
P_n(x)=\frac{1}{2}\sum_{j=0} ^{n+1}\left\{\binom{n+1}{j}x+\binom{n}{j-1}\right\}B_{n-1} ^{(n)}(x+j).
\end{equation*}
\end{theorem}
Let us consider the following two associated sequences:
\begin{equation}\label{53}
S_n(x)\sim\left(1,\frac{2t}{e^t+1}\right),~M_n(x)\sim\left(1,\frac{e^t-1}{e^t+1}\right).
\end{equation}
For $n \geq 1$, by \eqref{12}, we get
\begin{equation}\label{54}
\begin{split}
S_n(x)&=x\left(\frac{\frac{e^t-1}{e^t+1}}{\frac{2t}{e^t+1}}\right)^nx^{-1}M_n(x) \\
&=2^{-n}x\left(\frac{e^t-1}{t}\right)^nx^{-1}M_n(x).
\end{split}
\end{equation}
By \eqref{51} and \eqref{54}, we get
\begin{equation}\label{55}
\begin{split}
S_n(x)&=2^{-n}x\frac{1}{t^n}n!\sum_{l=n} ^{\infty}S_2(l,n)\frac{t^l}{l!}x^{-1}M_n(x) \\
&=2^{-n}x\sum_{l=0} ^{n-1} \frac{n!S_2(l+n,n)}{(l+n)!}t^l(x^{-1}M_n(x))\\
&=2^{-n}x\sum_{l=0} ^{n-1} \frac{n!}{(l+n)!}S_2(l+n,n)\sum_{j=0} ^n \binom{n}{j}(n-1)_lB_{n-1-l} ^{(n)} (x+j) \\
&=2^{-n}xn!(n-1)!\sum_{l=0} ^{n-1} \sum_{j=0} ^n \frac{\binom{n}{j}S_2(l+n,n)}{(l+n)!(n-l-1)!}B_{n-1-l} ^{(n)}(x+j).
\end{split}
\end{equation}
Therefore, by \eqref{50} and \eqref{53}, we obtain the following theorem.
\begin{theorem}\label{thm9}
For $n \geq 1$, we have
\begin{equation*}
\sum_{j=0} ^n \binom{n}{j}(x+j)^{n-1}=n!(n-1)!\sum_{l=0} ^{n-1} \sum_{j=0} ^n \frac{\binom{n}{j}S_2(l+n,n)}{(l+n)!(n-l-1)!}B_{n-1-l} ^{(n)}(x+j).
\end{equation*}
Moreover,
\begin{equation*}
\sum_{j=0} ^n \binom{n}{j}j^{n-1}=n!(n-1)!\sum_{l=0} ^{n-1} \sum_{j=0} ^n \frac{\binom{n}{j}S_2(l+n,n)}{(l+n)!(n-l-1)!}B_{n-1-l} ^{(n)}(j).
\end{equation*}
\end{theorem}
By \eqref{15}, we get
\begin{equation}\label{56}
x^{-1}M_n(x)=\sum_{k=1} ^n \binom{n}{k}(n-1)_{n-k}2^k\sum_{j=0} ^{k-1}S_1(k-1,j)(x-1)^j,
\end{equation}
where $S_1(k,j)$ is the Stirling number of the first kind. From \eqref{55} and \eqref{56}, we can derive
\begin{equation}\label{57}
\begin{split}
S_n(x)&=2^{-n}x\sum_{l=0} ^{n-1}\frac{n!}{(l+n)!}S_2(l+n,n)t^l(x^{-1}M_n(x))\\
&=2^{-n}x\sum_{l=0} ^{n-1}\frac{n!}{(l+n)!}S_2(l+n,n)t^l\sum_{k=1} ^n \binom{n}{k}(n-1)_{n-k}2^k \sum_{j=0} ^{k-1}S_1(k-1,j)(x-1)^j\\
&=2^{-n}xn!(n-1)!\sum_{j=0} ^{n-1}\sum_{l=0} ^{n-1} \sum_{k=j+1} ^n \frac{\binom{n}{k}2^kS_2(l+n,n)S_1(k-1,j)}{(l+n)!(k-1)!}t^l(x-1)^j\\
&=2^{-n}xn!(n-1)!\sum_{j=0} ^{n-1}\sum_{l=0} ^j\sum_{k=j+1} ^n \frac{\binom{n}{k}2^kS_2(l+n,n)S_1(k-1,j)}{(l+n)!(k-1)!}(j)_l(x-1)^{j-l}.
\end{split}
\end{equation}
Therefore, by \eqref{50} and \eqref{57}, we obtain the following theorem.
\begin{theorem}\label{thm10}
For $n \geq 1$, we have
\begin{equation*}
\sum_{j=0} ^n \binom{n}{j}(x+j)^{n-1}=n!(n-1)!\sum_{j=0} ^{n-1}\sum_{l=0} ^j\sum_{k=j+1} ^n \frac{\binom{n}{k}2^kS_2(l+n,n)S_1(k-1,j)j!}{(l+n)!(k-1)!(j-l)!}(x-1)^{j-l}.
\end{equation*}
\end{theorem}
{\scshape Remark.} From \eqref{51}, we note that
\begin{equation}\label{58}
x^{-1}M_n(x)=\sum_{k=0} ^n \binom{n}{k}(x+k-1)_{n-1}.
\end{equation}
By \eqref{55} and \eqref{58}, we get
\begin{equation}\label{59}
\begin{split}
S_n(x)&=2^{-n}x\sum_{l=0} ^{n-1}\frac{n!}{(l+n)!}S_2(l+n,n)t^l(x^{-1}M_n(x))\\
&=2^{-n}x\sum_{l=0} ^{n-1}\frac{n!}{(l+n)!}S_2(l+n,n)t^l\left(\sum_{k=0} ^n \binom{n}{k}(x+k-1)_{n-1}\right)\\
&=2^{-n}xn!\sum_{j=0} ^{n-1}\sum_{l=0} ^j \sum_{k=0} ^n \frac{\binom{n}{l}S_2(l+n,n)S_1(n-1,j)}{(l+n)!}(j)_l(x+k-1)^{j-l}\\
&=2^{-n}xn!\sum_{j=0} ^{n-1}\sum_{l=0} ^j\sum_{k=0} ^n \frac{\binom{n}{k}S_2(l+n,n)S_1(n-1,j)j!}{(l+n)!(j-l)!}(x+k-1)^{j-l}.
\end{split}
\end{equation}
So, by \eqref{50} and \eqref{59}, we get
\begin{equation}\label{60}
\sum_{j=0} ^n \binom{n}{j}(x+j)^{n-1}=n!\sum_{j=0} ^{n-1}\sum_{l=0} ^j\sum_{k=0} ^n \frac{\binom{n}{k}S_2(l+n,n)S_1(n-1,j)j!}{(l+n)!(j-l)!}(x+k-1)^{j-l}.
\end{equation}
The Narumi polynomials $N_n ^{(a)}(x)$ of order $a$ is defined by the generating function to be
\begin{equation}\label{61}
\sum_{k=0} ^{\infty}\frac{N_k ^{(a)}(x)}{k!}t^k=\left(\frac{\log (1+t)}{t}\right)^a(1+t)^x.
\end{equation}
Thus, from \eqref{61}, we see that
\begin{equation}\label{62}
N_n ^{(a)} (x)\sim\left(\left(\frac{e^t-1}{t}\right)^a,e^t-1\right),{\text{ (see [17,18])}}.
\end{equation}
In the special case, $x=0$, $N_k ^{(a)} (0)=N_k ^{(a)}$ are called the {\it{$k$-th Narumi numbers}} of order $a$. If $a=1$ in \eqref{62}, then we will write $N_n(x)$ and $N_n$ for $N_n ^{(1)} (x)$ and $N_n ^{(1)}$.

By \eqref{12} and \eqref{16}, we get
\begin{equation}\label{63}
\begin{split}
\phi_n(x)&=x\left(\frac{t}{\log (1+t)}\right)^nx^{-1}x^n=x\left(\frac{t}{\log (1+t)}\right)^nx^{n-1}\\
&=x\left(\sum_{k=0} ^{\infty} \frac{N_k ^{(-n)}(0)}{k!}t^k\right)x^{n-1}=x\sum_{k=0} ^{n-1} \frac{N_k ^{(-n)}}{k!}(n-1)_kx^{n-1-k}\\
&=\sum_{k=0} ^{n-1}\binom{n-1}{k}N_k ^{(-n)}x^{n-k}=\sum_{k=1} ^n \binom{n-1}{k-1} N_{n-k} ^{(-n)} x^k.
\end{split}
\end{equation}
Therefore, by \eqref{16} and \eqref{63}, we obtain the following lemma.
\begin{lemma}\label{lemma11}
For $n,k \in {\mathbb{N}}$ with $k \leq n$, we have
\begin{equation*}
S_2(n,k)=\binom{n-1}{k-1}N_{n-k} ^{(-n)}.
\end{equation*}
\end{lemma}
By \eqref{12}, \eqref{16} and \eqref{17}, we get
\begin{equation}\label{64}
\begin{split}
\phi_n(x)&=x\left(\frac{t}{(1+t)\log (1+t)}\right)^nx^{-1}L_n(-x) \\
&=x\left(\frac{t}{\log (1+t)}\right)^n(1+t)^{-n}\sum_{l=1} ^n \binom{n-1}{l-1}\frac{n!}{l!}x^{l-1} \\
&=x\left(\sum_{k=0} ^{\infty}\frac{N_k ^{(-n)} (-n)}{k!}t^k\right)\sum_{l=1} ^n \binom{n-1}{l-1}\frac{n!}{l!}x^{l-1} \\
&=x\sum_{k=0} ^{n-1}\frac{N_k ^{(-n)} (-n)}{k!}\sum_{l=k+1} ^n\binom{n-1}{l-1}\frac{n!}{l!}(l-1)_kx^{l-1-k}\\
&=n!\sum_{k=0} ^{n-1}\sum_{l=k+1} ^n\frac{\binom{n-1}{l-1}\binom{l-1}{k}}{l!}N_k ^{(-n)} (-n)x^{l-k}\\
&=n!\sum_{k=0} ^{n-1}\sum_{m=1} ^{n-k} \frac{\binom{n-1}{k+m-1}\binom{k+m-1}{k}}{(k+m)!}N_k ^{(-n)} (-n)x^m \\
&=n!\sum_{m=1} ^{n}\left\{\sum_{k=0} ^{n-m}\frac{\binom{n-1}{k+m-1}\binom{k+m-1}{k}}{(k+m)!}N_k ^{(-n)} (-n)\right\}x^m.
\end{split}
\end{equation}
From \eqref{16} and \eqref{64}, we have
\begin{equation}\label{65}
S_2(n,m)=n!\sum_{k=0} ^{n-m} \frac{\binom{n-1}{k+m-1}\binom{k+m-1}{k}}{(k+m)!}N_k ^{(-n)} (-n),
\end{equation}
where $1 \leq m \leq n$.

Therefore, by Lemma \ref{lemma11} and \eqref{65}, we obtain the following theorem.
\begin{theorem}\label{thm12}
For $m,n \in {\mathbb{N}}$ with $m \leq n$, we have
\begin{equation*}
\binom{n-1}{m-1}N_{n-m} ^{(-n)}=n!\sum_{k=0} ^{n-m}\frac{\binom{n-1}{k+m-1}\binom{k+m-1}{k}}{(k+m)!}N_k ^{(-n)} (-n).
\end{equation*}
\end{theorem}
It is well known that
\begin{equation}\label{66}
\left(\frac{t}{\log (1+t)}\right)^n(1+t)^{x-1}=\sum_{k=0} ^{\infty}B_k ^{(k-n+1)} (x)\frac{t^k}{k!},{\text{ (see [17])}}.
\end{equation}
Thus, by \eqref{61} and \eqref{66}, we get
\begin{equation}\label{67}
\sum_{k=0} ^{\infty}B_k ^{(k-n+1)} (x)\frac{t^k}{k!}=\left(\frac{t}{\log (1+t)}\right)^n(1+t)^{x-1}=\sum_{k=0} ^{\infty}N_k ^{(-n)}(x-1)\frac{t^k}{k!}.
\end{equation}
By comparing the coefficients on the both sides of \eqref{67}, we get
\begin{equation}\label{68}
B_k ^{(k-n+1)} (x)=N_k ^{(-n)} (x-1).
\end{equation}
Therefore, by \eqref{68}, we obtain the following corollary.
\begin{corollary}\label{coro13}
For $m,n\in {\mathbb{N}}$ with $m \leq n$, we have
\begin{equation*}
\binom{n-1}{m-1}B_{n-m} ^{(-m+1)} (1)=n!\sum_{k=0} ^{n-m} (\frac{{\binom{n-1}{k+m-1}\binom{k+m-1}{k}}}{{(k+m)!}}B_k ^{(k-n+1)}(-n+1).
\end{equation*}
\end{corollary}
Let us consider the following associated sequence:
\begin{equation}\label{69}
S_n(x)\sim\left(1,t(1+t)^a\right),~a\neq 0.
\end{equation}
Then, by \eqref{12} and \eqref{69}, we get
\begin{equation}\label{70}
\begin{split}
S_n(x)&=x\left(\frac{t}{t(1+t)^a}\right)^nx^{-1}x^n=x(1+t)^{-an}x^{n-1}\\
&=x\sum_{l=0} ^{n-1}\binom{-an}{l}t^lx^{n-1}=\sum_{l=0} ^{n-1} \binom{-an}{l}(n-1)_lx^{n-l} \\
&=\sum_{l=1} ^n \binom{-an}{n-l}(n-1)_{n-l}x^l.
\end{split}
\end{equation}
For $n \geq 1$, from \eqref{16} and \eqref{69}, we have
\begin{equation}\label{71}
\begin{split}
\phi_n(x)&=x\left(\frac{t(1+t)^a}{\log (1+t)}\right)^nx^{-1}S_n(x) \\
&=x\left(\frac{t}{\log (1+t)}\right)^n(1+t)^{an}x^{-1}S_n(x) \\
&=x\left(\sum_{k=0} ^{n-1} \frac{N_k ^{(-n)} (an)}{k!}t^k \right)x^{-1}S_n(x).
\end{split}
\end{equation}
By \eqref{70} and \eqref{71}, we get
\begin{equation}\label{72}
\begin{split}
\phi_n(x)&=x\sum_{k=0} ^{n-1}\frac{N_k ^{(-n)} (an)}{k!}\sum_{l=1} ^n \binom{-an}{n-l}(n-1)_{n-l}t^kx^{l-1} \\
&=x\sum_{k=0} ^{n-1}\frac{N_k ^{(-n)} (an)}{k!}\sum_{l=k+1} ^n \binom{-an}{n-l}(n-1)_{n-l}(l-1)_kx^{l-1-k} \\
&=(n-1)!\sum_{k=0} ^{n-1}\sum_{l=k+1} ^n\frac{\binom{-an}{n-l}\binom{l-1}{k}}{(l-1)!}N_k ^{(-n)}(an)x^{l-k} \\
&=(n-1)!\sum_{k=0} ^{n-1}\sum_{m=1} ^{n-k}\frac{\binom{-an}{n-k-m}\binom{k+m-1}{k}}{(k+m-1)!}N_k ^{(-n)}(an)x^m \\
&=(n-1)!\sum_{m=1} ^n\left\{\sum_{k=0} ^{n-m}\frac{\binom{-an}{n-k-m}\binom{k+m-1}{k}}{(k+m-1)!}N_k ^{(-n)}(an)\right\}x^m \\
&=(n-1)!\sum_{m=1} ^n\left\{\sum_{k=0} ^{n-m}\frac{\binom{-an}{n-k-m}\binom{k+m-1}{k}}{(k+m-1)!}B_k ^{(k-n+1)}(an+1)\right\}x^m.
\end{split}
\end{equation}
Therefore, by \eqref{16} and \eqref{72}, we obtain the following theorem.
\begin{theorem}\label{thm14}
For $m,n \in {\mathbb{N}}$ with $m \leq n$, we have
\begin{equation*}
\binom{n-1}{m-1}B_{n-m} ^{(-m+1)}(1)=(n-1)!\sum_{k=0} ^{n-m}\frac{\binom{-an}{n-k-m}\binom{k+m-1}{k}}{(k+m-1)!}B_k ^{(k-n+1)}(an+1).
\end{equation*}
\end{theorem}

{\scshape Remarks (I).} For $n \geq 1$, we have
\begin{equation}\label{73}
\begin{split}
x^n&=x\left(\frac{\log (1+t)}{t}\right)^nx^{-1}\phi_n(x) \\
&=\sum_{m=1} ^n \left\{\sum_{k=0} ^{n-m}\binom{k+m-1}{k}S_2(n,k+m)N_k ^{(n)}\right\}x^m.
\end{split}
\end{equation}
By comparing the coefficients on the both sides of \eqref{73}, we get
\begin{equation}\label{74}
\sum_{k=0} ^{n-m}\binom{k+m-1}{k}S_2(n,k+m)N_k ^{(n)}=\delta_{m,n},
\end{equation}
where $1 \leq m \leq n$.

{\scshape (II).} For $n \geq 1$, we have
\begin{equation}\label{75}
\begin{split}
L_n(-x)&=x\left(\frac{\log (1+t)}{\frac{t}{1+t}}\right)^nx^{-1}\phi_n(x) \\
&=\sum_{m=1} ^{n}\left\{\sum_{k=0} ^{n-m} \binom{k+m-1}{k}S_2(n,k+m)N_k ^{(n)} (n)\right\}x^m.
\end{split}
\end{equation}
By \eqref{17} and \eqref{75}, we get
\begin{equation}\label{76}
\binom{n-1}{m-1}\frac{n!}{m!}=\sum_{k=0} ^{n-m}\binom{k+m-1}{k}S_2(n,k+m)N_k ^{(n)} (n),
\end{equation}
where $1 \leq m \leq n$.

As is known, the Laguerre polynomials of order $\alpha$ are given by Sheffer sequences to be
\begin{equation}\label{77}
L_n ^{(\alpha)}(x)\sim\left((1-t)^{-\alpha-1},\frac{t}{t-1}\right).
\end{equation}
Thus, by the definition of Sheffer sequence, we get
\begin{equation}\label{78}
\left.\left<(1+t)^{-\alpha-1}\left(\frac{t}{t+1}\right)^k\right|L_n ^{(\alpha)}(-x)\right>=n!\delta_{n,k}~(n,k \geq 0).
\end{equation}
By \eqref{78}, we easily see that
\begin{equation}\label{79}
L_n(-x)=(1+t)^{-\alpha-1}L_n ^{(\alpha)}(-x)\sim\left(1,\frac{t}{1+t}\right).
\end{equation}
From \eqref{75} and \eqref{79}, we have
\begin{equation}\label{80}
\begin{split}
&(1+t)^{-\alpha-1}L_n ^{(\alpha)}(-x)=L_n(-x) \\
=&\sum_{m=1} ^{n}\left\{\sum_{k=0} ^{n-m} \binom{k+m-1}{k}S_2(n,k+m)N_k ^{(n)} (n)\right\}x^m.
\end{split}
\end{equation}
Thus, by \eqref{80}, we get
\begin{equation}\label{81}
\begin{split}
L_n ^{(\alpha)} (-x)&=(1+t)^{\alpha+1}\sum_{m=1} ^n \left\{\sum_{k=0} ^{n-m} \binom{k+m-1}{k}S_2(n,k+m)N_k ^{(n)} (n)\right\}x^m \\
&=\sum_{l=0} ^n \left\{\sum_{m=l} ^n \sum_{k=0} ^{n-m}\binom{k+m-1}{k}\binom{\alpha+1}{m-l}(m)_{m-l}S_2(n,k+m)N_k ^{(n)}(n)\right\}x^l.
\end{split}
\end{equation}
It is known that
\begin{equation}\label{82}
L_n ^{(\alpha)}(-x)=\sum_{l=0} ^n \binom{n+\alpha}{n-l}\frac{n!}{l!}x^l.
\end{equation}
By \eqref{81} and \eqref{82}, we get
\begin{equation*}
\binom{n+\alpha}{n-l}\frac{n!}{l!}=\sum_{m=l} ^n \sum_{k=0} ^{n-m} \binom{k+m-1}{k}\binom{\alpha+1}{m-l}(m)_{m-l}S_2(n,k+m)N_k ^{(n)}(n),
\end{equation*}
where $0\leq l \leq n$.

Finally, we consider the following associated sequences:
\begin{equation}\label{83}
S_n(x)=\sum_{k=1} ^n \binom{-an}{n-k}(n-1)_{n-k}x^k\sim\left(1,t(1+t)^a\right),~a\neq 0.
\end{equation}
Thus, by \eqref{12} and \eqref{83}, we get
\begin{equation}\label{84}
\begin{split}
S_n(x)&=x\left(\frac{\log (1+t)}{t(1+t)^a}\right)^nx^{-1}\phi_n(x) \\
&=x\left(\frac{\log (1+t)}{t}\right)^n(1+t)^{-an}x^{-1}\phi_n(x)\\
&=x\sum_{k=0} ^{n-1}\frac{N_k ^{(n)} (-an)}{k!}t^k\sum_{l=0} ^n S_2(n,l)x^{l-1} \\
&=\sum_{m=1} ^n \left\{\sum_{k=0} ^{n-m} \binom{k+m-1}{k}S_2(n,k+m)N_k ^{(n)} (-an)\right\}x^m.
\end{split}
\end{equation}
From \eqref{83} and \eqref{84}, we have
\begin{equation*}
\binom{-an}{n-m}(n-1)_{n-m}=\sum_{k=0} ^{n-m}  \binom{k+m-1}{k}S_2(n,k+m)N_k ^{(n)} (-an),
\end{equation*}
where $m,n \in {\mathbb{N}}$ with $m \leq n$ and $a \neq 0$.

\end{document}